% Logic Eprints
%Submitted 1725 Wed Dec 02, 1992 by: morass@math.mit.edu ()
%logic/friedman_sy/a_large_pi-1-2_set_absolute_for_set_forcing.tex
%

\input amstex
\documentstyle{amsppt}
\magnification=\magstep1
\NoBlackBoxes
\hsize 6.0truein
\vsize 8.2truein
\hcorrection{.2truein}
\vcorrection{.2truein}
\loadbold

\define\ORD{\operatorname{ORD}}
\define\Card{\operatorname{Card}}
\define\card{\operatorname{card}}

\document
\baselineskip=15pt
%\baselineskip=.25truein

\font\bigtenrm=cmr12 scaled\magstep2
\centerline
{\bigtenrm {A Large $\Pi^1_2$ Set, Absolute for Set Forcings}}

\vskip20pt

\font\bigtenrm=cmr10 scaled\magstep2
\centerline
{\bigtenrm{Sy D. Friedman}\footnote"*"{Research supported by NSF contract 
\#9205530-DMA.}} 

\centerline
{\bigtenrm {M.I.T.}}

\vskip20pt

\comment
lasteqno 1@0
\endcomment

The purpose of this note is to prove the following. 

\proclaim{Theorem}  Let $\kappa$ be an
$L$-cardinal, definable in $L.$  Then there is a set of reals $X,$
class-generic over $L,$  such that 

\vskip5pt

(a) \ $L(X)\models\Card=\Card^L$ and $X$ has cardinality $\kappa.$ 

(b) \ Some fixed $\Pi^1_2$ formula defines $X$  in all set-generic extensions
of $L(X).$ 
\endproclaim
\vskip5pt

\flushpar
By L\'evy-Shoenfield Absoluteness, any $\Pi^1_2$ formula defining $X$  in
$L(X)$ defines a superset of $X$  in each extension of $L(X).$  The point of
(b) is that this superset is just $X$  in set-generic extensions of $L(X).$
If $O^{\#}$ exists then $X$  as in the Theorem actually exists in $V,$  though
of course it will be only countable there.

The basic idea of the proof comes from David [82]. In his paper a real $R$
class-generic over $L$  is produced so that $\{R\}$ is $\Pi^1_2,$ uniformly
for set-generic extensions of $L(R).$  The added technique here is to use
``diagonal supports'' to take a large product of David-style forcings.

Here are some further applications of the Theorem and its proof.

\proclaim{Corollary 1} Assume consistency of an inaccessible cardinal. Then it
is consistent for the Perfect Set Property to hold for
$\underset\sim\to{\Sigma}^1_2$  sets yet fail for some $\Pi^1_2$ set.
\endproclaim

\demo{Proof} Using the Theorem get a $\Pi_2^1$ set $X$  which has cardinality
$\kappa$ in $L(X),$ $\kappa=$ least $L$-inaccessible, and which has a
$\Pi^1_2$-definition uniform for set-generic extensions. Then gently collapse
$\kappa$ to $\omega_1$ and add $\omega_2$ Cohen reals. In this extension,
$\omega_1>\omega_1^{L(R)}$ for each real $R$  and $X$  is a $\Pi^1_2$ set of
cardinality $\omega_1<\omega_2=2^{\aleph_0}.$  \hfill{$\dashv$ }
\enddemo

\proclaim{Corollary 2} Assume consistency of an inaccessible. Then it is
consistent that the Perfect Set Property holds for
$\underset\sim\to\Sigma^1_2$  sets and there is a $\Pi^1_2$ well ordering of
some set of reals of length $\aleph_{1000}.$  
\endproclaim

\vskip5pt

The latter answers a question of Harrington.

\newpage

\font\bigtenrm=cmr7 scaled\magstep2
\centerline
{\bigtenrm{The Proof.}} 

\vskip10pt

We modify the construction of David [82] to suit our purposes. First we
describe the $\alpha^+$-Souslin tree $T_\alpha$ in $L,$  where $\alpha$  is a
successor $L$-cardinal: $T_\alpha$  has a unique node on level $0$  and
exactly $2$  immediate successors on level $\beta+1$  to each node on level
$\beta,$  for $\beta<\alpha^+.$  If $\beta<\alpha^+$ is a limit of cofinality
$<\alpha$ then level $\beta$ assigns a top to each branch through the tree
below level $\beta.$  Now suppose $\beta<\alpha^+$ has cofinality $\alpha.$
Let ${\Cal{P}}$  be the forcing consisting of pairs $(\gamma,f)$ where
$\gamma<\beta$  and $f$  is a function from $\gamma$  into the nodes at levels
$<\beta,$  with extension defined by $(\gamma',f')\le(\gamma,f)$ iff
$\gamma'\ge\gamma,$  $f'(\delta)$ tree-extends $f(\delta)$ for each
$\delta<\gamma.$  Choose $G$  to be ${\Cal{P}}$-generic over $L_{\beta^*}$
where $\beta^*=$ largest p\.r\. closed $\beta^*\ge\beta$ such that
$\beta^*=\beta$ or $L_{\beta^*}\models\card(\beta)>\alpha.$ Then the nodes on
level $\beta$  are obtained by putting tops on the branches defined by
$\{f(\delta)|(\gamma,f)\in G$ some $\gamma\}$ for $\delta<\beta.$ 
This completes the definition of the $\alpha^+$-Souslin tree $T_\alpha.$

Now fix an $L$-definable cardinal $\kappa$ and also fix an $L$-definable $1-1$
function $F:\ \kappa\times\omega\times\ORD\longrightarrow$ Successor
$L$-cardinals greater than $\kappa.$  The forcing ${\Cal{P}}(\gamma,n),$
$\gamma<\kappa$ and $n<\omega,$ is designed to produce a real $R(\gamma,n)$
coding branches through $T_\alpha$  whenever $\alpha$  is of the form
$F(\gamma,n,\delta)$ for some $\delta.$  This  forcing is obtained by
modifying the Jensen coding of the empty class(see Beller-Jensen-Welch [82])
as follows: In defining the strings $s: [\alpha,|s|)\longrightarrow 2$ in
$S_\alpha,$  require that Even $(s)$ code a branch through $T_\alpha$ if
$\alpha\in\Card(\gamma,n)=\{F(\gamma,n,\delta)|\delta\in\ORD\}.$ Also use
David's trick to create a $\Pi^1_2$ condition implying that branches through
the appropriate trees are coded: for any $\alpha,$  for $s$  to belong to
$S_\alpha$  require that for $\xi\le|s|$ and $\eta>\xi,$ if
$L_\eta(s\restriction\xi)\models\xi=\alpha^++ZF^-+\Card=\Card^L$ then
$L_\eta(s\restriction\xi)\models$ for some $\gamma^*<\kappa^*,$
Even$(s\restriction\xi)$  codes a branch through $T^*_{\alpha^*}$  whenever
$\alpha^*\in\Card^*(\gamma^*,n),$ where $\kappa^*,$ $T^*_{\alpha^*},$
$\Card^*(\gamma^*,n)$ are defined in $L_\eta$  as were $\kappa,T_\alpha,$
$\Card(\gamma,n)$ in $L.$  The $\le\alpha$-distributivity of
${\Cal{P}}(\gamma,n)_\alpha(=$ the forcing at and above $\alpha)$  is
established as in David [82], with one added observation: if
$\alpha'\in\Card(\gamma,n)$ then we have to be sure that Even $(p_{\alpha'})$
codes a branch through $T_{\alpha'}$, where $p$  arises as the greatest lower
bound to an $\alpha$-sequence constructed to meet $\alpha$-many open dense
sets. There is no problem if $\alpha'>\alpha$ since then $T_{\alpha'}$ is
$\le\alpha$-closed. If $\alpha'=\alpha$ then the property follows from the
definition of level $|p_\alpha|$ of $T_\alpha,$ since we can arrange  that
Even $(p_\alpha)$ is sufficiently generic for $T_\alpha\restriction$ (levels
$<|p_\alpha|).$  (In fact the latter genericity is a consequence of the usual
construction of the $\alpha$-sequence leading to $p.)$ 

The forcing ${\Cal{P}}(\gamma), \gamma <\kappa,$ is designed to produce a real
$R(\gamma)$ such that $n\in R(\gamma)$ iff $R(\gamma)$ codes a branch
through $T_\alpha$ for each $\alpha$  in Card $(\gamma,n).$ A condition is
$p\in\prod\limits_{n}{\Cal{P}}(\gamma,n)$ where $p(n)(0)$ (a finite object) is
$(\emptyset\emptyset)$ for all but finitely many $n.$  Extension is defined by
$q\le p$ iff $q(n)\le p(n)$ in ${\Cal{P}}(\gamma,n)$ unless $n$ is not of the
form $2^{n_0}3^{n_{1}}$ or $n=2^{n_0}3^{n_1}$ where $q(n_0)_0(n_1)=0,$ in
which case there is no requirement on $q(n).$  A generic  $G$  can be
identified with the real $\{2^n3^m|p(n)_0(m)=1$ for some $p\in G\}=R(\gamma).$
The forcing at or above $\alpha,$ ${\Cal{P}}(\gamma)_\alpha,$ obeys
``quasi-distributivity'': if $D_i,i<\alpha$ are predense below $p$  then
there are  $q\le p$ and $d_i\subseteq D_i,i<\alpha$ such that each $d_i$ is
countable and predense below $q.$  This is established as in David [82] by
``guessing at $\langle p(n)(0)|n\in\omega\rangle$'' and yields cardinal
preservation. 

Our desired forcing ${\Cal{P}}$ is the ``diagonally supported'' product of the
${\Cal{P}}(\gamma),\gamma<\kappa.$  Specifically, a condition is
$p\in\prod\limits_{\gamma<\kappa }{\Cal{P}}(\gamma)$ where for infinite
cardinals $\alpha<\kappa,$
$\{\gamma|p(\gamma)(\alpha)\ne(\emptyset,\emptyset)\}$ has cardinality
$\le\alpha$ and in addition $\{\gamma|p(\gamma)(0)\ne(\emptyset,\emptyset)\}$
is finite. Quasi-distributivity for ${\Cal{P}}_\alpha=$ forcing at or above
$\alpha$  follows just as for ${\Cal{P}}(\gamma)_\alpha.$ The point of the
diagonal supports is that for infinite successor cardinals $\alpha,{\Cal{P}}$
factors as ${\Cal{P}}_\alpha*{\Cal{P}}^{G_{\alpha }}$ where $G_\alpha$ denotes
the ${\Cal{P}}_\alpha$-generic and ${\Cal{P}}^{G_{\alpha }}$ is $\alpha^+-CC.$
Thus we get cardinal-preservation.

Now note that if $\langle R(\gamma)|\gamma<\kappa\rangle$ comes from  (and
therefore determines) a ${\Cal{P}}$-generic then $n\in
R(\gamma)\longrightarrow R(\gamma)$ codes a branch through $T_\alpha$ for
$\alpha$  in Card $(\gamma,n).$ Conversely, if $n\notin R(\gamma)$ then there
is no condition on extension of conditions in ${\Cal{P}}(\gamma)$ to cause
$R(\gamma)$ to code a branch through such $T_\alpha.$  In fact, by the
quasi-distributivity argument for ${\Cal{P}}_\alpha,$  given any term
$\tau$ for a subset of $\alpha^+$ and any condition $p,$  we can find
$\beta<\alpha^+$ of cofinality $\alpha$  and $q\le p$ such that $q$ forces
$\tau\cap\beta$ to be one of $\alpha$-many possibilities, each constructed
before $\beta^*,$ where $\beta=|q_\alpha|.$  Thus $q$ forces that $\tau$ is
{\it not} a branch through $T_\alpha,$ so we get: $n\in R(\gamma)$  iff
$R(\gamma)$ codes a branch through each $T_\alpha,$
$\alpha\in\Card(\gamma,n)$ iff $R(\gamma)$ codes a branch through some
$T_\alpha,$ $\alpha\in\Card(\gamma,n).$  The coding is localized in the sense
that if $n\in R(\gamma)$ then whenever $L_\eta(R(\gamma))\models
ZF^-+\Card=\Card^L,$ there is $\gamma^*<\kappa^*$ such that
$L_\eta(R(\gamma))\models R(\gamma)$ codes a branch through $T^*_{\alpha^*}$
whenever $\alpha^*\in\Card^*(\gamma^*,n),$ where
$\kappa^*,T^*_{\alpha^*},\Card^*(\gamma^*,n)$ are defined in $L_\eta$ just as
$\kappa,T_\alpha,$ $\Card(\gamma^*,n)$ are defined in $L.$  The latter
condition on $R(\gamma)$ is sufficient to know that $R(\gamma)$ is equal to
one of the intended $R(\gamma),\gamma<\kappa,$ even if we restrict ourselves
to countable $\eta.$  With that restriction we get a $\Pi^1_2$ condition
equivalent to membership in $X=\{R(\gamma)|\gamma<\kappa\}.$ Since set-forcing
preserves the Souslin-ness of trees at sufficiently large cardinals, the above
$\Pi^1_2$ definition of $X$ works in any set-generic extension of $L(X).$
This completes the proof of the Theorem.

\demo{Proof of Corollary 2} As in the proof of Corollary 1 we can obtain
$X=\{R(\gamma)|\gamma<\kappa\},$ $\kappa=999^{\text{th}}$ cardinal after the
least 
$L$-inaccessible, which has a $\Pi^1_2$ definition uniform for set-generic
extensions of $L(X),$ where $\Card^{L(X)}=\Card^L.$ We can guarantee that
$Y=\{\langle R(0), R(\gamma_1),
R(\gamma_2)\rangle|0<\gamma_1\le\gamma_2<\kappa\}$ also has such a uniform
$\Pi^1_2$ definition, using the following trick: Design $R(0)$ so that $u\in
R(0)\Longleftrightarrow$  Even $(R(0))$ codes a branch through $T_\alpha$ for
each $\alpha$ in Card $(0,n),$ and so that Odd$(R_0)$ almost disjointly codes
$\{\langle R(\gamma_1),R(\gamma_2)\rangle|0<\gamma_1\le\gamma_2<\kappa\}.$
Thus, for $R\in L(X),R^*$ is almost disjoint from Odd$(R_0)$ iff $R=\langle
R(\gamma_1),R(\gamma_2)\rangle$  for some $0<\gamma_1\le\gamma_2<\kappa,$
where $R^*=\{n|n$ codes a finite initial segment of $R\}.$  The former
requires only a very small modification to the definition of the
${\Cal{P}}(0)$ forcings. The latter requires only a small modification to the
definition of ${\Cal{P}}:$ take the diagonally-supported product as before,
but restrain $p(0)$ for $p\in{\Cal{P}}$ so as to affect the desired almost
disjoint coding. These finite restraints do not interfere with the
quasi-distributivity argument for ${\Cal{P}}.$ 

Now we have the desired $\Pi^1_2$ definition for $Y=\{\langle R(0),
R(\gamma_1), R(\gamma_2)\rangle|0<\gamma_1\le\gamma_2<\kappa\}:R$ belongs to
$Y$ iff $R=\langle R_0, R_1, R_2\rangle$  where $R_0=R(0)$  and $\langle
R_1,R_2\rangle^*$ is almost disjoint from $R_0$ and $R_1,R_2$ belong to $X.$
Since $R(0)$ is uniformly definable as a $\Pi^1_2$-singleton in set-generic
extensions of $L(X),$ this is the desired definition. Of course, using $Y$  we
obtain a $\Pi^1_2$ well-ordering of length $\kappa.$  Finally as in the proof
of Corollary 1, gently collapse $\kappa$ to $\omega_1$ and we have
$\omega_1>\omega_1^{L(R)}$ for each real $R$ with a $\Pi^1_2$ well-ordering of
length $\aleph_{1000}.$ \hfill{$\dashv$ }
\enddemo

\vskip10pt

\flushpar
{\bf Remarks.} \ The same proof gives length $\aleph_\alpha$ for any
$L$-definable $\alpha.$  We can also add Cohen reals so that the continuum is
as large as desired, without changing the maximum length of a $\Pi^1_2$
well-ordering.

It is possible to show that if $O^{\#}$ exists then there is  a $\Pi^1_2$ set
$X$  such that $X$  has large cardinality in $L(X).$  But this requires the
more difficult technique of Friedman [90].

\enddocument